\documentclass[10pt]{amsart}
\usepackage{amsthm, amsfonts, amssymb, amsmath, amscd}
\usepackage[all]{xy}
\usepackage{graphicx}
\usepackage{lmodern}
\usepackage[margin=3.5cm]{geometry}

\newcommand{\To}{\rightarrow}

\newcommand{\into}{\hookrightarrow}

\newcommand{\noin}{\noindent}

\newcommand{\C}{\mathbb{C}}
\newcommand{\R}{\mathbb{R}}

\newcommand{\Z}{\mathbb{Z}}

\newcommand{\Cstar}{\mathrm{C^*}}
\newcommand{\Cred}{\mathrm{C^*_r}}
\newcommand{\G}{\Gamma}

\newcommand{\KH}{\mathcal{K}}
\newcommand{\U}{\mathrm{U}}
\newcommand{\M}{\mathrm{M}}
\newcommand{\GL}{\mathrm{GL}}
\newcommand{\Idem}{\mathrm{Idem}}
\newcommand{\any}{\:\cdot\:}

\newcommand{\tsr}{\mathrm{tsr}\:}

\newcommand{\bsr}{\mathrm{bsr}\:}
\newcommand{\csr}{\mathrm{csr}\:}

\newcommand{\gsr}{\mathrm{gsr}\:}

\newcommand{\spec}{\mathrm{sp}}

\newcommand{\dom}{\mathrm{dom}}

\newcommand{\bij}{\mathrm{bij}}

\newcommand{\id}{\mathrm{id}}

\newcommand{\Lip}{\mathrm{Lip}}
\newcommand{\Lg}{\mathrm{Lg}}

\newtheorem{thm}{Theorem}[section]

\newtheorem{prop}[thm]{Proposition}

\theoremstyle{definition}

\newtheorem{ex}[thm]{Example}

\newtheorem*{ack}{Acknowledgments}

\begin{document}
\title{Spectral morphisms, K-theory, and stable ranks}
\author{Bogdan Nica}
\address{\newline Department of Mathematics and Statistics, University of Victoria, Victoria (BC), Canada V8W 3R4}
\date{\today}
\dedicatory{To the memory of Israel M. Gelfand (1913 - 2009)}

\begin{abstract}
\noin We give a brief account of the interplay between spectral morphisms, K-theory, and stable ranks in the context of Banach algebras.
\end{abstract}

\thanks{Supported by a Postdoctoral Fellowship from the Pacific Institute for the Mathematical Sciences (PIMS)}
\maketitle

\section{Introduction}
This article revolves around the notion of spectral morphism - that is, a morphism which preserves spectra of elements - in the setting of Banach algebras. The importance of spectral morphisms was firmly established, in the commutative case, by the Gelfand transform. In the general, non-commutative case, spectral morphisms with dense image have been very useful for K-theoretic purposes: the Karoubi - Swan theorem asserts that topological K-theory is invariant across such morphisms. Since K-theory is controlled by stable ranks - which are, very loosely speaking, noncommutative notions of dimension - the following question arises: are stable ranks invariant across spectral morphisms with dense image? In the literature, this is known as Swan's problem. Schematically, the inter-connections between spectral morphisms, K-theory, and stable ranks are summarized by the following diagram:

\vspace{1cm}
\begin{equation*}
\xymatrix{
&\ar@/_1.2pc/[ldd]^{\genfrac{}{}{0pt}{}{\textrm{Karoubi - Swan}}{\textrm{theorem}}} & \genfrac{}{}{0pt}{}{\textrm{\large spectral}}{\textrm{\large morphisms}} &\ar@/^1.2pc/[rdd]_{\genfrac{}{}{0pt}{}{\textrm{Swan's}}{\textrm{problem}}}  &\\
& &  \ar@{.>}@/^1.2pc/[ldd]^{\genfrac{}{}{0pt}{}{\textrm{Gelfand}}{\textrm{transform}}} \ar@{.>}@/_1.2pc/[rdd] & &\\
 && & &\\
 \genfrac{}{}{0pt}{}{\textrm{\large topological}}{\textrm{\large K-theory}} &  && &\genfrac{}{}{0pt}{}{\textrm{\large \;\; stable \;\;}}{\textrm{\large \;\; ranks \;\;}}\\
&   &  & \ar@/^1.4pc/[ll]_{\genfrac{}{}{0pt}{}{\textrm{stabilization for homotopy}}{\textrm{groups of $\GL_n$}}} &}
\end{equation*}

\vspace{1cm}

At the Noncommutative Geometry Workshop (May 2008), I gave a talk based on my paper \cite{Nic08}. A straightforward summary of the paper seemed uninteresting. So I decided to interpret the opportunity offered by these Proceedings as an occasion to describe the broader context captured in the above diagram. Results from \cite{Nic08}, and from a more recent paper \cite{Nic09}, are interwoven in this picture.

We follow the interplay between spectral morphisms and K-theory, on one hand, and spectral morphisms and stable ranks, on the other hand, in two different contexts. The first context is the Gelfand transform. The second context, typically non-commutative, deals with spectral morphisms with dense image; informally, the density assumption is supposed to make up for the loss of commutativity. We shall witness a nearly perfect parallelism between these two contexts.

The last section is devoted to stabilization phenomena for the homotopy groups of the general linear group over Banach algebras. The stable ranks play the key role in this discussion.

\section{Spectral morphisms}
In what follows, Banach algebras and their (continuous) morphisms are assumed to be unital. 

A Banach algebra morphism $\phi: A\to B$ is \emph{spectral} if it is spectrum-preserving, that is, for all $a\in A$ we have $\spec_B (\phi(a))=\spec_A (a)$. Equivalently, the morphism $\phi$ is spectral if, for all $a\in A$, we have that $a$ is invertible in $A$ if and only if $\phi(a)$ is invertible in $B$. 

In the case when $\phi$ is an inclusion - which is the typical case - $A$ is spectral in $B$ if and only if $A$ is holomorphically closed in $B$; see, for instance, \cite{Sch92}. Recall, $A$ is said to be holomorphically closed in $B$ if, for each $a\in A$, we have that $f(a)\in A$ whenever $f$ is holomorphic in a neighborhood of $\spec_B(a)$.  The property that $A$ is spectral in $B$ appears under many other different names in the literature: Bourbaki (\cite{Bou67}) calls $A$ full in $B$; Naimark (\cite{Nai72}) says that $A$ and $B$ form a Wiener pair; other sources say that $A$ is inverse-closed in $B$.

The following example of spectral morphism is fundamental:
\begin{ex}[Gelfand transform]
Let $A$ be a unital, commutative Banach algebra. The \emph{maximal ideal space} $X_A$ is the set of characters of $A$, equipped with the topology of pointwise convergence; the name is due to the fact that we have a bijective correspondence, given by $\chi\mapsto\ker\chi$, between the characters of $A$ and the maximal ideals of $A$. The maximal ideal space $X_A$ is a compact Hausdorff space. 

Let $\hat{a}\in C(X_A)$ denote the evaluation at $a\in A$. Then the evaluation map $a\mapsto \hat{a}$ defines a unital, continuous morphism $A\To C(X_A)$. This morphism - the \emph{Gelfand transform} of $A$ - is a spectral morphism. 
\end{ex}

Outside of the Gelfand context, the spectral morphisms we consider are assumed to have dense image. Henceforth, a morphism with dense image is called, simply, a \emph{dense} morphism. The Gelfand transform may or may not be dense; an instance when it is dense is when the image is $*$-closed.

If $A$ and $B$ are commutative and $A\to B$ is a spectral morphism, then a simple argument involving the determinant shows that $\M_n(A)\to\M_n(B)$ is spectral for each $n\geq 1$. A noncommutative analogue is provided by the following useful fact, due to Swan \cite{Swa77} (see also Bost \cite{Bos90} and Schweitzer \cite{Sch92}):

\begin{prop}\label{matrixspectral}
If $A\to B$ is a dense and spectral morphism, then $\M_n(A)\to\M_n(B)$ is a dense and spectral morphism for each $n\geq 1$.
\end{prop}

Let us give some examples of morphisms which are dense and spectral. Especially motivating from the point of view of noncommutative geometry is the following commutative example:

\begin{ex}\label{smooth}
Let $M$ be a compact smooth manifold. Then $C^k(M)$ is a Banach algebra under the norm $\|f\|_{(k)}:=\max\{\|\partial^\alpha f\|_\infty: 0\leq |\alpha|\leq k\}$ (defined using local charts on $M$), and the inclusion $C^k(M) \into C(M)$ is dense and spectral.
\end{ex}

Passing from the smooth category to the metric category, we still have an analogue of the Banach algebra of $C^1$-functions on a compact manifold:

\begin{ex} Let $X$ be a compact metric space. Then the algebra $\Lip (X)$ of Lipschitz functions on $X$ is a Banach algebra under the norm $f\mapsto \max\{\|f\|_\infty, \|f\|_{\Lip}\}$, and the inclusion $\Lip (X)\into C(X)$ is dense and spectral.
\end{ex}

Let us point out that the dense and spectral inclusions appearing in the previous examples are Gelfand transforms. This is due to the fact that a dense and spectral morphism $A\to B$ between commutative Banach algebras induces a homeomorphism $X_B\to X_A$; in particular, a dense and spectral morphism $A\to C(X)$, where $X$ is a compact Hausdorff space, can be identified with the Gelfand transform of $A$.

The next two examples provide dense and spectral Banach subalgebras for the reduced $\Cstar$-algebras of certain groups. Recall that the reduced $\Cstar$-algebra $\Cred\G$ of a group $\G$ is the completion of the group algebra $\C\G$ under the operator norm  $\|\cdot\|$ coming from the regular representation of $\G$ on $\ell^2\G$.

\begin{ex}\label{PolyGrowth}
Let $\G$ be a finitely generated group. If $\G$ has polynomial growth, then the dense inclusion $\ell^1\G\into\Cred\G$ is spectral (Ludwig \cite{Lud79}).
\end{ex}

\begin{ex}\label{RD}
Let $\G$ be a finitely generated group, equippped with a word-length $|\cdot|$. Following Jolissaint \cite{Jol90}, we define the $s$-Sobolev space $H^s\G$ to be the Banach space obtained by completing $\C\G$ under the weighted $\ell^2$-norm 
\[\big\|\sum a_gg\big\|_{2,s}:=\sqrt{\sum |a_g|^2(1+|g|)^{2s}}.\] 
The group $\G$ is said to have \emph{property RD (of order $s$)} if, for some $C,s\geq 0$, we have $\|a\|\leq C\|a\|_{2,s}$ for all $a\in \C\G$. Free groups (Haagerup \cite{Haa79}) - more generally, hyperbolic groups (de la Harpe \cite{dHa88}) - have property RD; so do groups of polynomial growth (Jolissaint \cite{Jol90}), and many other groups.

If $\G$ has property RD of order $s$, then the following holds: for every $S>s$, the $S$-Sobolev space $H^S\G$ is in fact a Banach algebra under $\|\cdot\|_{2,S}$, and the (continuous) dense inclusion $H^S\G\into \Cred\G$ is spectral (Lafforgue \cite{Laf00}).
\end{ex}

The following example, of a more general nature, is the noncommutative analogue of Example~\ref{smooth}:

\begin{ex} \label{derivations}
Let $B$ be a Banach algebra and $\delta: B\To B$ an unbounded derivation, i.e., $\delta$ is a linear map defined on a subalgebra $\dom (\delta)$ of $B$ and satisfying the Leibniz rule $\delta(ab)=\delta(a)b+a\delta(b)$. Also, assume that $\delta$ is closed; then, for each $k\geq 1$, $B^{(k)}:=\dom (\delta^k)$ is a Banach algebra under the norm $\|a\|_{(k)}:=\max\{\|\delta^i(a)\|: 0\leq i\leq k\}$. 

If $B^{(k)}$ is dense in $B$, then the (continuous) inclusion $B^{(k)}\into B$ is spectral (Ji \cite{Ji92}).
\end{ex}

As a non-example, we mention the following counterpoint to Example~\ref{PolyGrowth}:

\begin{ex}\label{FS2} Let $\G$ be a group containing a free subsemigroup on two generators. Then the dense inclusion $\ell^1\G\into\Cred\G$ is not spectral (Jenkins \cite{Jen70}).
\end{ex}

\section{Spectral morphisms and K-theory}\label{SM-K}
We recall the definition of the K-groups, the main purpose being that of introducing some notation. Following Taylor \cite{Tay76}, we use the shorthand $[\any]$ for the set of path-connected components $\pi_0(\any)$.

Let $A$ be a unital Banach algebra. Let $V(A): =\varinjlim\; \big[\Idem(\M_n(A))\big]$, where $\Idem(\any)$ denotes the set of idempotents and the direct limit is taken along the maps induced by $*\mapsto \big(\begin{smallmatrix} *&0\\0&0 \end{smallmatrix}\big)$. Then $V(A)$ is an abelian monoid under $[a]+[a']:=\big[ \big(\begin{smallmatrix} a&0\\0&a' \end{smallmatrix}\big) \big]$. Also, let $\mathcal{P}(A)$ denote the set of isomorphism classes of finitely generated, projective right $A$-modules; equipped with the direct sum, $\mathcal{P}(A)$ is an abelian monoid. The monoids $V(A)$ and $\mathcal{P}(A)$ are naturally isomorphic, and the Grothendieck group of $V(A)$ - or $\mathcal{P}(A)$ - is the (abelian) group $K_0(A)$. On the other hand, we simply put $K_1(A):=\varinjlim\; \big[\GL_n(A)\big]$, where $\GL_n(A)$ denotes the invertible group of $\M_n(A)$ and the direct limit is taken along the maps induced by $*\mapsto \big(\begin{smallmatrix} *&0\\0&1 \end{smallmatrix}\big)$. This time, the direct limit is already an abelian group under $[a][a']:=\big[ \big(\begin{smallmatrix} a&0\\0&a' \end{smallmatrix}\big) \big]$.

\subsection{The Gelfand transform} Consider the Gelfand transform $A\To C(X_A)$ of a commutative Banach algebra $A$. Historically, the first result with a K-theoretic flavor is Shilov's idempotent theorem \cite{Shi53} saying that the induced map $\Idem(A)\To\Idem(C(X_A))$ is onto. It is easy to see, on the other hand, that the Gelfand transform is injective on idempotents: if $p,q$ are idempotents with $\hat{p}=\hat{q}$, then $z:=pq+(1-p)(1-q)$ is invertible since $\hat{z}=1$, and from $pz=zq(=pq)$ we deduce - by the commutativity of $A$ - that $p=q$. Thus, an equivalent form of Shilov's theorem is the statement that the Gelfand transform $A\To C(X_A)$ induces a bijection $\Idem(A)\To\Idem(C(X_A))$. In hindsight, this gives a strong hint that a $K_0$-isomorphism might be lurking in the background. 

The next K-theoretic development, this time in the direction of $K_1$, is the theorem of Arens \cite{Are63} and Royden \cite{Roy63}: the Gelfand transform $A\To C(X_A)$ induces an isomorphism $\big[A^\times\big]\To \big[C(X_A)^\times\big]$, where $(\any)^\times$ denotes the invertible group. Shortly after, Arens \cite{Are66} extended the Arens - Royden theorem by showing that, for each $n\geq 1$, we have an isomorphism $\big[\GL_n(A)\big]\To \big[\GL_n(C(X_A))\big]$ induced by the Gelfand transform $A\To C(X_A)$.

Clearly, the Arens theorem implies that the Gelfand transform $A\To C(X_A)$ induces an isomorphism $K_1(A)\to K_1(C(X_A))$. This was first pointed out by Novodvorskii \cite{Nov67}, who also handled the $K_0$-isomorphism for semisimple $A$. The semisimple hypothesis was later removed by Taylor; one of the results from \cite{Tay76} is that the Gelfand transform $A\To C(X_A)$ induces, for each $n\geq 1$, a bijection $\big[\Idem(\M_n(A))\big]\To\big[\Idem(\M_n(C(X_A)))\big]$. Summarizing, we have 

\begin{thm}[Novodvorskii - Taylor]\label{Kisom - Gelfand}
The Gelfand transform $ A\To C(X_A)$ induces an isomorphism $K_*(A)\to K_*(C(X_A))$.
\end{thm}

The Arens theorem has been greatly generalized by Davie \cite{Dav71}. As a particular case of Davie's theorem, we have the following fact: the Gelfand transform $A\To C(X_A)$ induces a homotopy equivalence $\GL_n(A)\To \GL_n(C(X_A))$ for each $n\geq 1$.

\subsection{Dense and spectral morphisms} 
Versions of the following theorem are due to Karoubi \cite{Kar78} and Swan \cite{Swa77}; see also Connes \cite{Con81}. An excellent account is given by Bost in \cite[Appendix]{Bos90}. 

\begin{thm}[Karoubi - Swan]\label{Kisom - SM}
A dense and spectral morphism $A\To B$ induces an isomorphism $K_*(A)\to K_*(B)$.
\end{thm}

The Karoubi - Swan theorem features, for instance, in the proof by Connes and Moscovici \cite{CM90} that hyperbolic groups satisfy the Novikov conjecture, and in the remarkable work of Lafforgue \cite{Laf02}. In both cases, one needs to pass from the reduced $\Cstar$-algebra $\Cred \G$ of a group $\G$ to a suitable group algebra of $\G$ in such a way that the K-theory is preserved.

The proof of the Karoubi - Swan theorem is fairly elementary. One shows, in fact, that a dense and spectral morphism $A\To B$ induces bijections $\big[\GL_n(A)\big]\To \big[\GL_n(B)\big]$ and $\big[\Idem(\M_n(A))\big]\To \big[\Idem(\M_n(B))\big]$ for all $n\geq 1$. In light of Proposition~\ref{matrixspectral}, it suffices to consider the case $n=1$. Now, to prove that $\big[\Idem(A)\big]\To \big[\Idem(B)\big]$ is a bijection, one first enlarges the closed set of idempotents $\Idem (\any)$ to a suitable open set of ``almost'' idempotents; for instance, one could take the set $^\sim \Idem(\any)$ consisting of all elements whose spectrum does not meet the line $\mathrm{Re}\; z=\frac{1}{2}$. By holomorphic functional calculus, the set of idempotents $\Idem (\any)$ is a deformation retract of the set of almost idempotents $\mathrm{^\sim Idem}(\any)$; therefore, it suffices to prove that $\big[\mathrm{^\sim Idem}(A)\big]\To \big[\mathrm{^\sim Idem}(B)\big]$ is a bijection. This, and the fact that $\big[A^\times\big]\To \big[B^\times\big]$ is a bijection, are proved in a similar manner. 

Actually, these two bijections are avatars of the following ``spectral'' picture (cf. \cite{Nic08}). Let $\Omega$ be an open subset of $\C$ containing the origin. For a unital Banach algebra $A$, we let $A_\Omega$ denote the set of elements of $A$ whose spectrum is contained in $\Omega$. For example, $A_\Omega=A^\times -1$ when $\Omega=\C\setminus \{-1\}$, and $A_\Omega=\mathrm{^\sim Idem}(A)$ when $\Omega=\{z:\; \mathrm{Re}\; z\neq\frac{1}{2}\}$. Then a dense and spectral morphism $A\To B$ induces a bijection $\big[A_\Omega\big]\To \big[B_\Omega\big]$. As soon as we adopt this spectral picture, the K-theoretic pull is to define a ``spectral'' functor $K_\Omega$ in the same way as the usual K-functors are defined. Given a unital Banach algebra $A$, the direct limit
\[V_\Omega(A): =\varinjlim\; \big[\M_n(A)_\Omega\big]\]
is an abelian monoid under $[a]+[a']=\big[ \big(\begin{smallmatrix} a&0\\0&a' \end{smallmatrix}\big) \big]$. Applying the Grothendieck construction to the abelian monoid $V_\Omega(A)$, we obtain an abelian group $K_\Omega (A)$. Thus, for each open set $\Omega\subseteq \C$ containing the origin, we have a \emph{spectral K-functor} $K_\Omega$. Of course, not every choice of $\Omega$ yields an interesting functor; but, at the very least, we can recover the classical functors $K_0$ and $K_1$. The Karoubi - Swan theorem for spectral K-functors reads as follows: a dense and spectral morphism $A\To B$ induces an isomorphism $K_\Omega (A)\to K_\Omega (B)$. The interesting part is not this generalized Karoubi - Swan theorem, but rather how the Karoubi - Swan theorem lead us to the generalized K-functors $K_\Omega$. 

The analogue, due to Bost (\cite{Bos90}), of Davie's result is the expected one: a dense and spectral morphism $A\To B$ induces, for each $n\geq 1$, a homotopy equivalence $\GL_n(A)\To \GL_n(B)$. In fact, with the previous notations we have that such a $\phi$ induces, for each $n\geq 1$, a homotopy equivalence $\M_n(A)_\Omega\To \M_n(B)_\Omega$.

\section{Spectral morphisms and stable ranks}
Bass \cite{Bas64} introduced the first notion of stable rank; more than a dozen stable ranks have appeared since then. Like the original stable rank introduced by Bass, some are purely ring-theoretic; other stable ranks take topology into account. Here, we consider the following four stable ranks: the Bass stable rank, the topological stable rank, the connected stable rank, and the general stable rank. The last three stable ranks are due to Rieffel \cite{Rie83}.

Let $A$ be a unital Banach algebra. For $n\geq 1$, consider the set of \emph{unimodular $n$-tuples}
\[\Lg_n(A)=\{(a_1,\dots,a_n):\; Aa_1+\dots + Aa_n=A\}\subseteq A^n.\]
The last column of an invertible $n$-by-$n$ matrix is an example of unimodular $n$-tuple. Left-multiplying the transpose of a unimodular $n$-tuple by an invertible $n$-by-$n$ matrix yields another 
unimodular $n$-tuple; this describes an action of $\GL_n(A)$ on $\Lg_n(A)$. Under this action, the orbit of $(0,\dots,0,1)$ consists of the last columns of invertible matrices.

The first two of the four stable ranks under consideration are defined as follows. The \emph{Bass stable rank} of $A$, denoted $\bsr A$, is the least $n\geq 1$ with the property that $\Lg_{n+1}(A)$ is ``reducible'' to $\Lg_n(A)$ in the following sense: if $(a_1,\dots, a_{n+1})\in \Lg_{n+1}(A)$, then $(a_1+x_1a_{n+1}, \dots, a_n+x_na_{n+1})\in \Lg_n(A)$ for some $(x_1, \dots, x_n)\in A^n$. The \emph{topological stable rank} of $A$, denoted $\tsr A$, is the least $n\geq 1$ with the property that $\Lg_n(A)$ is dense in $A^n$. Naturally, a stable rank is declared to be infinite if there is no integer $n$ satisfying the required condition. The Bass and the topological stable ranks are related by the inequality $\bsr A\leq \tsr A$ (\cite{Rie83}, \cite{CL84}).

If $X$ is a compact space then 
\[\bsr C(X)=\tsr C(X)=\Big\lfloor \frac{\dim X}{2} \Big\rfloor +1\]
where $\dim X$ denotes the covering dimension of $X$. In light of this formula - due to Vaserstein \cite{Vas71} for the Bass stable rank, and to Rieffel \cite{Rie83} for the topological stable rank - we think of the Bass and the topological stable ranks as being \emph{dimensional} stable ranks.

The remaining two stable ranks are defined as follows. The \emph{connected stable rank} of $A$, denoted $\csr A$, is the least $n\geq 1$ with the property that $\GL^0_m(A)$ acts transitively on $\Lg_m(A)$ for all $m\geq n$. By $\GL^0_m(A)$ we denote, as usual, the identity component of $\GL_m(A)$. It is not hard to see that the action of $\GL^0_m(A)$ on $\Lg_m(A)$ has open orbits. Hence, the \emph{connected stable rank} of $A$ can also be defined as the least $n\geq 1$ with the property that $\Lg_m(A)$ is connected for all $m\geq n$. Finally, the \emph{general stable rank} of $A$, denoted $\gsr A$, is the least $n\geq 1$ with the property that $\GL_m(A)$ acts transitively on $\Lg_m(A)$ for all $m\geq n$. It turns out (Rieffel \cite{Rie83}) that an equivalent definition is the following: the \emph{general stable rank} of $A$ is the least  $n\geq 1$ such that, for all $m\geq n$, if $P$ is a right $A$-module satisfying $P\oplus A\simeq A^m$ then $P\simeq A^{m-1}$. Clearly, we have that $\gsr A\leq \csr A$.

The connected and the general stable ranks are invariant under homotopy equivalence (\cite{Nis86}, \cite{Nic09}): if two Banach algebras $A$ and $B$ are homotopy equivalent (i.e., there are Banach algebra morphisms $\phi:A\to B$ and $\psi: B\to A$ such that $\psi\phi$ is homotopic to $\id_A$ and $\phi\psi$ is homotopic to $\id_B$ via paths of Banach algebra morphisms) then $\csr A=\csr B$ and $\gsr A=\gsr B$. For this reason, we think of the  connected and the general stable ranks as being \emph{homotopical} stable ranks (cf. \cite{Nic09}). The similarities between these two stable ranks, and the contrast between the homotopical and the dimensional stable ranks is the main theme of \cite{Nic09}. One such difference is that the homotopical stable ranks of $C(X)$ can be quite hard to compute, in general. Furthermore, here we encounter a difference amongst the homotopical stable ranks: $\csr C(X)$ can be computed whenever $X$ satisfies a fairly mild cohomological hypothesis (\cite{Nic08}, \cite{Nic09}), but no such general result is known for $\gsr C(X)$. In fact, even the computation of $\gsr C(T^d)$ - where $T^d$ denotes the familiar $d$-dimensional torus - is not known!

Stable ranks turn out to be relevant in K-theory; see Section~\ref{sr-K}. As we have seen in Section~\ref{SM-K}, K-theory is preserved by suitable spectral morphisms. The following question arises: are stable ranks also preserved by these spectral morphisms? This question, first suggested in \cite{Swa77}, is called \emph{Swan's problem}.

\subsection{The Gelfand transform} We start with the behavior of the dimensional stable ranks under the Gelfand transform $A\To C(X_A)$ of a commutative Banach algebra $A$. For the Bass stable rank, we have $\bsr A\leq \bsr C(X_A)$ (Corach - Larotonda \cite{CL84}); examples such as the disk algebra $A(D)$ or the Hardy algebra $H^\infty(D)$, where $D$ denotes the open unit disk, show that the inequality can be strict. However, if the Gelfand transform $A\To C(X_A)$ is also dense then $\bsr A=\bsr C(X_A)$ (Vaserstein \cite{Vas71}).  For the topological stable rank, no general comparison between $\tsr A$ and $\tsr C(X_A)$ is known.

On the other hand, the homotopical stable ranks are invariant across the Gelfand transform:

\begin{thm}[Novodvorskii - Taylor, Forster - Taylor]\label{Gelfand csrgsr} Let $A\To C(X_A)$ be the Gelfand transform of a commutative Banach algebra $A$. Then $\csr A=\csr C(X_A)$ and $\gsr A=\gsr C(X_A)$.
\end{thm}

The invariance of the connected stable rank is due to the fact that the Gelfand transform $A\To C(X_A)$ induces a bijection $\big[\Lg_n(A)\big]\to \big[\Lg_n(C(X_A))\big]$ for all $n\geq 1$. When $A$ is semisimple, this fact is a consequence of Novodvorskii's results from \cite{Nov67}; in general, it is proved by Taylor in \cite{Tay76}.

The first result pointing towards the invariance of the general stable rank is due to Lin \cite{Lin73} (see also Taylor \cite{Tay76}) who showed that $(a_1,\dots, a_n)\in \Lg_n(A)$ is the last column of a matrix in $\GL_n(A)$ if and only if $(\hat{a}_1,\dots, \hat{a}_n)\in \Lg_n(C(X_A))$ is the last column of a matrix in $\GL_n (C(X_A))$; in particular, we have that $\gsr A\leq\gsr C(X_A)$. To obtain the desired invariance of the general stable rank, recall instead the argument which showed the invariance of $K_0$ across the Gelfand transform: $A\To C(X_A)$ induces a bijection $\big[\Idem(\M_n(A))\big]\To\big[\Idem(\M_n(C(X_A)))\big]$ for all $n\geq 1$, hence a monoid isomorphism $\mathcal{P}(A)\to\mathcal{P}(C(X_A))$ (Forster \cite{For74}, Taylor \cite{Tay76}). Since the general stable rank can also be defined in terms of a certain cancellation property for projective modules, it trivially follows that the general stable rank is invariant across the Gelfand transform. 

\subsection{Dense and spectral morphisms} We start again with the dimensional stable ranks: the Bass stable rank and the topological stable rank. If $A\To B$ is a dense and spectral morphism, then 
\[\bsr A\leq \bsr B \leq \tsr B\leq \tsr A.\] 
The inequality $\tsr B\leq \tsr A$ is easily proved under the mere assumption that $A\To B$ has dense image. Since spectrality plays no role, this inequality cannot be taken as a hint for what should happen in the Gelfand context. 

The inequality $\bsr A\leq \bsr B$, due to Swan \cite{Swa77}, mirrors (and predates) the Corach - Larotonda result from the Gelfand context. An instance when equality holds is provided by the following result of Badea \cite{Bad98}: if $A$ is a dense and spectral Banach $*$-subalgebra of a $\Cstar$-algebra $B$, and $A$ is closed under $C^\infty$-functional calculus for self-adjoint elements, then $\bsr A=\bsr B$; this applies, for example, to the dense and spectral subalgebras constructed from derivations (Example~\ref{derivations}). 

Having $\bsr A= \tsr A$ would, of course, solve Swan's problem for both dimensional stable ranks. This is known to hold whenever $A$ is a $\Cstar$-algebra (Herman - Vaserstein \cite{HV}); unfortunately, the domain of a dense and spectral morphism is a $\Cstar$-algebra only in trivial cases. Therefore, extensions of the Herman - Vaserstein theorem outside of the scope of $\Cstar$-algebras would be very relevant here.

For the homotopical stable ranks, namely the connected stable rank and the general stable rank, we have a very satisfactory result:

\begin{thm}\label{SM-csrgsr}
If $A\To B$ is a dense and spectral morphism, then $\csr A=\csr B$ and $\gsr A=\gsr B$.
\end{thm}

The invariance of the connected stable rank was first proved in \cite{Nic08}, as a consequence of the fact that a dense and spectral morphism $A\To B$ induces a bijection $\big[\Lg_n(A)\big]\to \big[\Lg_n(B)\big]$ for all $n\geq 1$; in \cite{Nic09}, a different proof is given. The invariance of the general stable rank can be proved by a direct argument (see \cite{Nic09}), or by invoking the following analogue of the Forster - Taylor theorem from the Gelfand context: a dense and spectral morphism $A\To B$ induces a monoid isomorphism $\mathcal{P}(A)\to\mathcal{P}(B)$; this follows from the fact that $\big[\Idem(\M_n(A))\big]\To \big[\Idem(\M_n(B))\big]$ is a bijection for all $n\geq 1$. However, the direct argument has the added benefit of being adaptable to the ring-theoretic context of Swan \cite{Swa77}.

\section{Relatively spectral morphisms}
Starting from the notion of spectral morphism, there are several ways of running into a weaker property. Firstly, consider a dense and spectral morphism of Banach algebras $A\to B$.  Full surjectivity is not needed - on the $B$ side, one approximates by elements from the image of $A$. This idea, that density suffices, could then be applied on the $A$ side: the spectral condition could be required to hold on a dense subalgebra of $A$ only. Secondly, instead of a dense inclusion consider a nesting of dense inclusions of Banach algebras $A\into B\into C$. As $\spec_A(a)\supseteq \spec_B(a)\supseteq \spec_C(a)$ for all $a\in A$, the inclusion $A\into C$ is spectral if and only if the inclusion $A\into B$ is spectral and the inclusion $B\into C$ is spectrum-preserving over the dense subalgebra $A$. Thirdly, let $\G$ be a group and consider the dense inclusion $\ell^1\G\into\Cred\G$, or any other dense inclusion involving two Banach algebra completions of the group algebra $\C\G$. Amidst the hazy halo of such a Banach algebra completion, the elements of $\C\G$ have a reassuring concreteness. One is compelled to think of the possibility that the spectrum-invariance property for the dense inclusion $\ell^1\G\into\Cred\G$ is only known for elements of $\C\G$, instead of having it satisfied throughout $\ell^1\G$.

These situations lead us to the following definition (cf. \cite{Nic08}). A Banach algebra morphism $\phi: A\to B$ is \emph{relatively spectral} if it is spectrum-preserving over a dense subalgebra of $A$. That is, for some dense subalgebra $X$ of $A$ we have $\spec_B (\phi(x))=\spec_A (x)$ for all $x\in X$; equivalently, for all $x\in X$ we have that $x$ is invertible in $A$ if and only if $\phi(x)$ is invertible in $B$. 

One naturally wonders at this point:
\begin{itemize}
\item[($\dagger$)] Let $A\to B$ be a relatively spectral morphism of Banach algebras. Is $A\to B$ spectral?
\end{itemize}
The answer to ($\dagger$) is positive whenever $A$ satisfies some form of spectral continuity (see \cite{Nic08}). It is unlikely, however, that the property in ($\dagger$) is always true. As a particular aspect of ($\dagger$), we do not know whether the relative analogue of Proposition~\ref{matrixspectral} holds. So let us call a Banach algebra morphism $A\to B$ \emph{completely relatively spectral} if each amplified morphism $\M_n(A)\to \M_n(B)$ is relatively spectral. This (slightly awkward) terminology follows the standard practice of calling a morphism $A\to B$ \emph{completely $\mathcal{P}$} whenever the amplified morphisms $\M_n(A)\to \M_n(B)$ satisfy a certain property $\mathcal{P}$.

Here are the relative analogues of Examples~\ref{PolyGrowth} and~\ref{derivations}:

\begin{ex}\label{SubexpGrowth}
Let $\G$ be a finitely generated group. If $\G$ has subexponential growth, then the dense inclusion $\ell^1\G\into\Cred\G$ is completely relatively spectral.
\end{ex}

\begin{ex} \label{derivations, in general}
Let $B$ be a Banach algebra and $\delta: B\To B$ an unbounded derivation. For $k\geq 1$, we let $B^{(k)}$ denote the Banach algebra obtained by completing the (unital) algebra $\dom (\delta^k)$ under the algebra norm $\|a\|_{(k)}:=\max\{\|\delta^i(a)\|: 0\leq i\leq k\}$. 

If $B^{(k)}$ is dense in $B$, then the (continuous) inclusion $B^{(k)}\into B$ is completely relatively spectral.
\end{ex}

It is worth pointing out that Example~\ref{SubexpGrowth} is much easier to prove than Example~\ref{PolyGrowth}. Concerning Example~\ref{SubexpGrowth}, we also mention that the question whether $\ell^1\G\into\Cred\G$ is spectral whenever $\G$ has subexponential growth - an instance of ($\dagger$) - is equivalent to the question whether $\ell^1\G$ is symmetric whenever $\G$ has subexponential growth. The latter question, still not settled, is a part of the following old problem: characterize the groups $\G$ for which the Banach algebra $\ell^1\G$ is symmetric. Recall, a unital Banach $*$-algebra is \emph{symmetric} if every self-adjoint element has real spectrum.

As for Example~\ref{derivations, in general}, it is a generalization of Example~\ref{derivations}. Indeed, one shows that each inclusion $\M_n(B^{(k)})\into \M_n(B)$ is spectral relative to the dense subalgebra $\M_n(\dom (\delta^k))$; if $\delta$ is assumed to be closed, then $\dom (\delta^k)$ coincides with $B^{(k)}$ and we recover the result of Example~\ref{derivations}.

Why care about relatively spectral morphisms? The point is that we can reap the same benefits as in the case of spectral morphisms, but with less spectral information. Namely, a Banach algebra morphism which is dense and completely relatively spectral induces an isomorphism in K-theory and preserves the homotopical stable ranks (\cite{Nic08}, \cite{Nic09}).

\section{Stable ranks and K-theory}\label{sr-K}
\subsection{Stabilization in algebraic K-theory} Before we get to stabilization in topological K-theory, let us provide some motivation coming from the algebraic side. Bass \cite{Bas64} devised his ``stable range'' condition as a way to control the sequence 
\[\dots\to\GL_n(R)/\mathrm{E}_n(R)\to\GL_{n+1}(R)/\mathrm{E}_{n+1}(R)\to\cdots\]
which yields the algebraic $K_1$-group $K^\textrm{alg}_1(R)$ in the limit. Recall that, for a unital ring $R$, the elementary group $\mathrm{E}_n(R)$ is the subgroup of $\GL_n(R)$ generated by the elementary matrices $\{1_n+e_{ij}\}_{1\leq i\neq j\leq n}$. We think of the quotients $\GL_n(R)/\mathrm{E}_n(R)$ as being the non-stable $K_1$-groups of $R$, despite the fact that $\GL_n(R)/\mathrm{E}_n(R)$ is not necessarily a group - that is, $\mathrm{E}_n(R)$ may not be normal in $\GL_n(R)$. The map $\GL_n(R)\to\GL_{n+1}(R)$, given by $*\mapsto \big(\begin{smallmatrix} *&0\\0&1 \end{smallmatrix}\big)$, sends $\mathrm{E}_n(R)$ to $\mathrm{E}_{n+1}(R)$, so it induces a well-defined map $\GL_n(R)/\mathrm{E}_n(R)\to\GL_{n+1}(R)/\mathrm{E}_{n+1}(R)$.

The following theorem was conjectured by Bass in \cite{Bas64}. Bass proved the surjectivity part, whereas the injectivity part, significantly more difficult, is due to Vaserstein \cite{Vas69}.

\begin{thm}[Bass - Vaserstein] 
The map $\GL_n(R)/\mathrm{E}_n(R)\to\GL_{n+1}(R)/\mathrm{E}_{n+1}(R)$
is surjective for $n\geq \bsr R$, and injective for $n\geq \bsr R+1$.
\end{thm}
In particular, $\GL_n(R)/\mathrm{E}_n(R)$ is a group isomorphic to $K^\mathrm{alg}_1(R)$ for $n\geq \bsr R+1$.

In \cite{Sus}, Suslin extended the Bass - Vaserstein theorem to higher algebraic K-theory. For a sequence of non-stable $K_i$-groups
\[K_{i,1}(R)\to K_{i,2}(R)\to\dots\to K_{i,n}(R)\to K_{i,n+1}(R)\to\cdots\]
we expect stabilization whenever $n$ is large enough compared to $\bsr R$ and $i$. As it is usually the case with higher algebraic K-theory, a significant problem is how to define these non-stable $K_i$-groups. The two possibilities considered by Suslin are $K^V_{i,n}(R)$ in the sense of Volodin, and $K^Q_{i,n}(R)$ in the sense of Quillen; see \cite{Sus} for the definitions. There is a canonical map $K^V_{i,n}(R)\to K^Q_{i,n}(R)$ for $n\geq 2i+1$. For these two interpretations of the non-stable $K_i$-groups, Suslin proved the following stabilization results:
\begin{enumerate}
\item[($K_i^V$)] The map $K^V_{i,n}(R)\to K^V_{i,n+1}(R)$ is surjective for $n\geq \bsr R + i-1$, and injective for $n\geq \bsr R+i$.

\item[($K_i^Q$)] The map $K^Q_{i,n}(R)\to K^Q_{i,n+1}(R)$ is surjective for $n\geq \max\{i+1, \bsr R\} + i-1$, and injective for $n\geq \max\{i+1, \bsr R\} + i$.
\end{enumerate}
A weaker version of ($K_i^Q$) was also obtained by Maazen and van der Kallen \cite{vdK80}.

\subsection{Stabilization in topological K-theory} Let $A$ be a (unital) Banach algebra. For each $k\geq 0$, we have a direct sequence of identity-based homotopy groups:
\[(\pi_k)\qquad\pi_k(\GL_1(A))\to\pi_k(\GL_2(A))\to\dots\to \pi_k(\GL_n(A))\to\pi_k(\GL_{n+1}(A))\to\cdots\]
The maps in this sequence are induced by the embeddings $\GL_n(A)\into\GL_{n+1}(A)$ given by $*\mapsto \big(\begin{smallmatrix} *&0\\0&1 \end{smallmatrix}\big)$. The direct limit of the homotopy sequence $(\pi_k)$ is, by definition, the topological K-group $K_{k+1}(A)$. By Bott periodicity, $K_{k+1}(A)$ is isomorphic to $K_0(A)$ or $K_1(A)$ according to whether $k$ is odd or even. 

We are interested in the following stabilization problem: when do the maps in the homotopy sequence $(\pi_k)$ become - and remain - injective / surjective / bijective? The first stabilization result in this direction appears in \cite{Rie83}. We formulate Rieffel's result in a way which makes it the topological analogue of the Bass - Vaserstein theorem:

\begin{thm}[Rieffel]\label{Rieffel I}
The map $\pi_0(\GL_{n}(A))\To\pi_0(\GL_{n+1}(A))$ is surjective for $n\geq \bsr A$, and injective for $n\geq \bsr A+1$.
\end{thm}

Actually, Rieffel shows that surjectivity holds for $n\geq \csr A-1$; since $\bsr A\geq\csr A-1$ (\cite{Rie83}, \cite{CL84}), surjectivity also holds for $n\geq \bsr A$. As for the fact that injectivity holds for $n\geq \bsr A+1$, Rieffel proves this by invoking the injectivity part of the Bass - Vaserstein theorem. 

Corach - Larotonda \cite{CL86} then generalized Theorem~\ref{Rieffel I} to higher homotopy groups as follows:

\begin{thm}[Corach - Larotonda]\label{CorachLarotonda} 
The map $\pi_k(\GL_{n}(A))\To\pi_k(\GL_{n+1}(A))$ is surjective for $n\geq \bsr A+k$, and injective for $n\geq \bsr A+k+1$.
\end{thm}

The proof of this result holds the key to all stabilization results in the literature. It relies on the existence of a fibration $\GL_n(A)\to \GL_{n+1}(A)\to\Lg_{n+1}(A)$, generalizing the more familiar fibration $\U_n\to \U_{n+1}\to S^{2n+1}$. For $n\geq \csr A-1$, $\Lg_{n+1}(A)$ is connected  and we get a long exact  homotopy sequence 
\[\dots \to\pi_{k+1}(\Lg_{n+1}(A))\to\pi_k(\GL_{n}(A))\to\pi_k(\GL_{n+1}(A))\to\pi_k(\Lg_{n+1}(A))\to\cdots\]
ending in
\[\dots\to\pi_1(\Lg_{n+1}(A))\to\pi_0(\GL_{n}(A))\to\pi_0(\GL_{n+1}(A))\to 0.\]
We identify $\pi_k(\Lg_{n+1}(A))$ with $\pi_0(\Lg_{n+1}(A(S^k)))$, where $A(S^k)$ denotes the Banach algebra of continuous maps from the $k$-sphere $S^k$ to $A$. This brings the connected stable rank of $A(S^k)$ into play, for which - as one can easily show - we have $\csr A(S^k)\geq \csr A$.
From the long exact homotopy sequence we conclude:

\begin{thm}\label{highercsr}
The map $\pi_k(\GL_{n}(A))\To\pi_k(\GL_{n+1}(A))$ is surjective for $n\geq \csr A(S^k)-1$, and injective for $n\geq \csr A(S^{k+1})-1$.
\end{thm}

To obtain Theorem~\ref{CorachLarotonda} from Theorem~\ref{highercsr}, one needs to know that $\bsr A+k+1\geq\csr A(S^k)$. To that end, it suffices to know that $\bsr A+k\geq \bsr A(S^k)$, and this is precisely what Corach - Larotonda show in the remainder of their proof.

In the case $k=0$, Theorem~\ref{highercsr} together with the fact that $\tsr A+1\geq \csr A(S^1)$ (Rieffel \cite{Rie83}) yield the following significant result from \cite{Rie87}:

\begin{thm}[Rieffel]\label{Rieffel II}
The map $\pi_0(\GL_{n}(A))\To\pi_0(\GL_{n+1}(A))$ is bijective for $n\geq \tsr A$.
\end{thm}

Comparing the two theorems of Rieffel, Theorem~\ref{Rieffel I} and Theorem~\ref{Rieffel II}, we see that the latter strengthens the former when $A$ is a $\Cstar$-algebra - recall, the Bass and the topological stable ranks coincide in this case. I do not know whether this is still the case for general Banach algebras.

A variant of the Corach - Larotonda fibration, where the base space $\Lg_n(A)$ is replaced by the subspace $\mathrm{Lc}_n(A)$ consisting of last columns of invertible $n$-by-$n$ matrices, is used by Rieffel in  \cite{Rie87} and later on by Thomsen \cite{Tho91} and Badea \cite{Bad99}, among others. This allows for the following sharpening of the injectivity bound from Theorem~\ref{highercsr}: the map $\pi_k(\GL_{n}(A))\To\pi_k(\GL_{n+1}(A))$ is injective for $n\geq \gsr A(S^{k+1})-1$ (Badea \cite{Bad99}). A downside of this improvement is that the general stable rank is, in general, more difficult to compute than the connected stable rank.

Compared to the Corach - Larotonda theorem, Theorem~\ref{highercsr} has a clear formal advantage. On the one hand, the injective / surjective stabilization of $(\pi_k)$ has a three-fold invariance: with respect to homotopy equivalence, across Gelfand transforms, and across dense and spectral morphisms. On the other hand, the bounds from Theorem~\ref{highercsr} enjoy the same three-fold invariance. A concrete example highlighting the power of Theorem~\ref{highercsr} over the Corach - Larotonda theorem is provided by the irrational rotation $\Cstar$-algebra $A_\theta$ (cf. \cite{Rie87}). It turns out that $\csr A_\theta (S^k)\leq 2$ for each $k\geq 0$; by Theorem~\ref{highercsr}, for each $k\geq 0$ we have that $\pi_k(\GL_n(A_\theta))$ is isomorphic to $K_{k+1}(A_\theta)\simeq \Z^2$ for all $n\geq 1$. Such a result would not be visible in its entirety through the Corach - Larotonda theorem.

All the results mentioned so far indicate an index $n$ from which all maps in the homotopy sequence $(\pi_k)$ are isomorphisms. But how to compute the smallest index $n$ with the property that the sequence $(\pi_k)$ stabilizes from $n$ onwards? The non-triviality of such a computation is already suggested by the trivial case $A=\C$: it can be shown that the smallest such index is $n=\lfloor k/2\rfloor+1$, but one needs rather detailed information about the homotopy groups of unitary groups for this.

Let $\bij_k\:A$ denote the least $n\geq 1$ such that the map $\pi_k(\GL_m(A))\To\pi_k(\GL_{m+1}(A))$ is bijective for all $m\geq n$; equivalently, $\bij_k\:A$ is the least $n\geq 1$ such that the natural map $\pi_k(\GL_m(A))\To K_{k+1}(A)$ is bijective for all $m\geq n$. By Theorem ~\ref{highercsr}, we have the following upper bound:
\[\bij_k\: A\leq \max\{\csr A(S^k)-1, \csr A(S^{k+1})-1\}\]
But the long exact homotopy sequence arising from the Corach - Larotonda fibration also yields the following inequality:
\[\csr A(S^{k+1})-1\leq \max\{\bij_k\:A,\bij_{k+1}\: A,\csr A-1\}\]
A variant of this inequality, formulated in terms of certain higher connected stable ranks, appears in \cite{Nic08}. These two inequalities allow us to compute $\max\{\bij_k\:A,\bij_{k+1}\: A\}$, which encodes the \emph{joint} stabilization of $(\pi_k)$ and $(\pi_{k+1})$. This stabilization for pairs of consecutive homotopy sequences has the following K-theoretic interpretation: if we think of stabilization in a single homotopy sequence $(\pi_k)$ as stabilization for one of the two K-groups, then the joint stabilization of $(\pi_k)$ and $(\pi_{k+1})$ should be understood as stabilization for the K-theory, i.e. for both K-groups, of $A$.

We first consider the case of commutative $\Cstar$-algebras. Let $A=C(X^d)$, where $X^d$ is a $d$-dimensional compact space, and assume that $H^d(X^d;\Z)\neq 0$. If $k\geq 0$ has the same parity as $d$, then $n=1+\frac{1}{2}(d+k)$ is the least integer with the property that the natural maps
\[\pi_{k}(\GL_m(A)) \to K_{k+1}(A),\qquad \pi_{k+1}(\GL_m(A)) \to K_{k+2}(A)\]
are isomorphisms for all $m\geq n$.

As a non-commutative example, we consider tensor products of $\Cstar$-extensions of compacts by commutative $\Cstar$-algebras. For $1\leq i\leq N$, let $A_i$ be a unital $\Cstar$-extension of $\KH$ by $C(X_i)$, where $X_i$ is a compact metrizable space. Put $A=A_1\otimes\dots \otimes A_N$, and $X=X_1\times\dots\times X_N$. Assume that $d=\dim\: X$ is non-zero, and $H^d(X^d; \Z)\neq 0$. If $k\geq 0$ has the same parity as $d$, then $n=1+\frac{1}{2}(d+k)$ is the least integer with the property that the natural maps
\[\pi_{k}(\GL_m(A)) \to K_{k+1}(A),\qquad \pi_{k+1}(\GL_m(A)) \to K_{k+2}(A)\]
are isomorphisms for all $m\geq n$.

Roughly speaking, this means that $A$ has the same stabilization pattern as its ``symbol space'' $C(X)$. The key point in both examples is, of course, the  computation of $\csr A(S^k)$. This is immediate for the first example, but the computation for the second example is more involved. The ingredients of this computation, performed in \cite{Nic09}, are estimates of stable ranks for quotients and extensions, as well as Nistor's result from \cite{Nis86} on the dimensional stable rank for such tensor products.

\begin{ack} I thank the referee for some useful suggestions.
\end{ack}

\end{document}